\documentclass[11pt]{amsart}
\usepackage{amssymb,latexsym,comment,url,amsmath}

\usepackage[T1]{fontenc}
\usepackage{tikz}
\usepackage{rotating,graphicx}
\usepackage{currvita}
\usepackage{diagbox} 
\usepackage{xcolor}
\usepackage{subcaption}
\usepackage{multirow}
\usepackage{array}
\usepackage{longtable}
\usepackage{ltablex}
\usepackage{float}
\usepackage{lmodern,babel,adjustbox,booktabs,multirow}
\usepackage[T1]{fontenc}
\usepackage[utf8]{inputenc}
\usepackage{ upgreek }
\usepackage{caption}
\usepackage{placeins}
\restylefloat{table}
\usepackage{blindtext}
\usepackage{hyperref}

\newcommand{\End}{{\operatorname{End}}}
\newcommand{\Aut}{{\operatorname{Aut}}}

\newcommand{\C}{{\mathbb{C}}}

\newcommand{\Char}{\operatorname{char}}

\newcommand{\Jac}{\operatorname{Jac}}

\newcommand{\Gal}{\operatorname{Gal}}

\newcommand{\isom}{ \cong }

\newcommand{\PP}{{\mathbb P}}
\newcommand{\Q}{{\mathbb Q}}

\newcommand{\Z}{{\mathbb Z}}

\newenvironment{Proof}{\par\noindent{\sc Proof:}}%
                      {\hspace*{\fill}\nobreak$\Box$\par\medskip}
                       {\hspace*{\fill}\nobreak$\Box$\par\medskip}

\newtheorem{Proposition}{Proposition}[section]
\newtheorem{Theorem}[Proposition]{Theorem}
\newtheorem{Lemma}[Proposition]{Lemma}

\newtheorem{Corollary}[Proposition]{Corollary}

\makeglossary
\newtheorem*{theorem*}{Theorem}

\theoremstyle{definition}
\newtheorem{Definition}[Proposition]{Definition}
\newtheorem{Remark}[Proposition]{Remark}
 \newtheorem{Example}[Proposition]{Example}

\renewcommand{\baselinestretch}{1.1}

\begin{document}

\title[Explicit construction of decomposable Jacobians]%
{Explicit construction of decomposable Jacobians}

\author[M. Buğday]%
{Mesut~Buğday}
\address{Faculty of Engineering and Natural Sciences, Sabanc{\i} University, Tuzla, \.{I}stanbul, 34956 Turkey}
\email{mesut.bugday@sabanciuniv.edu}
\author[M. Sadek]%
{Mohammad~Sadek}
\email{mohammad.sadek@sabanciuniv.edu}

\begin{abstract}
In this note we give explicit constructions of decomposable  hyperelliptic Jacobian varieties over fields of characteristic $0$. These include hyperelliptic Jacobian varieties that are isogenous to a product of two absolutely simple hyperelliptic Jacobian varieties, a square of a hyperelliptic Jacobian variety, and a product of four hyperelliptic Jacobian varieties three of which are of the same dimension. As an application, we produce families of hyperelliptic curves with infinitely many quadratic twists having at least two rational non-Weierstrass points; and families of quadruples of hyperelliptic curves together with infinitely many square-free $d$ such that the quadratic
twists of each of the curves by $d$ possess at least one rational non-Weierstrass point.
\end{abstract}
\maketitle

\let\thefootnote\relax

\section{Introduction}
An abelian variety is said to be {\em decomposable} over a field $K$ if it is isogenous to a product of abelian varieties of lower dimension. The study of decomposable Jacobian varieties of genus two curves was initiated in \cite{Haya}, see also \cite{Kuhn}.  
A family of hyperelliptic curves of arbitrary genus whose Jacobians decompose into two abelian varieties was given in \cite{Earle}, namely,
for the Jacobian of the hyperelliptic curve defined by the equation
$$y^2 = (x^n-1)(x^n-t), \quad n=2k+1,\quad k>1,\quad t\in\C\setminus \{0,1\},$$
there are two algebraic curves $Y_1$ and $Y_2$ of genus $k$ such that $\Jac(X)$ is isomorphic to $\Jac(Y_1) \times \Jac(Y_2)$. Ekedahl
and Serre constructed examples of curves whose Jacobians decompose completely into elliptic curves, \cite{Ekedahl}. The reader may also see \cite{Yam} for such examples of curves over number fields. Jacobian varieties of algebraic curves with many automorphisms provide examples of abelian varieties that contain many factors in their decompositions. In \cite{Paulhos, Jennifer, Jennifer2}, such curves whose Jacobians contain many elliptic factors were displayed. In \cite{Carocca}, the existence of Jacobians that are isogenous to
the product of arbitrary many Jacobians of the same genus, not necessarily equal to one, was established. 

In this note, we consider the following question. Given a positive integer $n$ together with a partition $n_1\le n_2\le ...\le n_k$ of $n$, does there exist a Jacobian variety of dimension $n$ that decomposes into a product of $k$ Jacobian varieties of dimensions $n_1,\cdots, n_k$? When $k=2$ and $n$ is even, we give explicit examples of families of hyperelliptic Jacobian varieties that decompose into the product of two absolutely simple Jacobian varieties of the same dimension $n/2$; and families of hyperelliptic Jacobian varieties that decompose as the square of a Jacobian variety. When $n$ is odd, we present examples of hyperelliptic Jacobian varieties that decompose into the product of two absolutely simple Jacobian varieties of dimensions $(n-1)/2$ and $(n+1)/2$. We exhibit families of hyperelliptic Jacobians that decompose into the product of three Jacobians of dimensions $k$, $k+1$, $2k$ when $n=4k+1$, $k\ge 1$; and $k+1$, $k+1$, $2k+1$ when $n=4k+3$, $k\ge0$. Further, we prove the existence of hyperelliptic Jacobian varieties of odd dimension $n$ that decompose as the product of four Jacobian varieties of dimensions $k$, $k$, $k$, $k+1$, when $n=4k+1$, $k\ge 1$; and $k$, $k+1$, $k+1$, $k+1$ when $n=4k+3$, $k\ge1$. In particular, given any integer $M$, there is a decomposable Jacobian variety of dimension $4M\pm 1$ whose decomposition contains three Jacobian factors each of dimension $M$. 

Goldfeld Conjecture states that the average rank of elliptic curves over the rational field in families
of quadratic twists is $1/2$. In other words,
quadratic twists of an elliptic curve over the rational field with rank at least $2$ are rare. In \cite{Rubin, Kuwata}, quadratic twists of elliptic curves with ranks at least $2$ or $3$ were given. A similar problem was posed to find tuples of elliptic curves whose quadratic twists by the same rationals are of positive rank infinitely often, \cite{Alaa, Coogan, Im}. As for hyperelliptic curves, one may construct families of these curves with infinitely many quadratic twists that possess no rational points, \cite{Sadek1,Sadek2,Legrand}. As a byproduct of our construction of decomposable Jacobian varieties, we produce examples of hyperelliptic curves with infinitely many quadratic twists possessing at least two rational non-Weierstrass points. In particular, we introduce examples of elliptic curves with infinitely many quadratic twists of rank at least $2$. In addition, we give examples of families of quadruples of hyperelliptic curves, three of which are of the same genus, such that for infinitely many square-free rationals the quadratic twists of each of these hyperelliptic curves by these rationals possess at least one rational non-Weierstrass point.

\subsection*{Acknowledgment}
The authors would like to thank the anonymous referee for many suggestions that improved the manuscript. These suggestions include the statement and proof of Theorem \ref{thm:ref}.

This work is supported by The Scientific and Technological Research Council of Turkey, T\"{U}B\.{I}TAK, research grant ARDEB 1001/122F312. M. Sadek acknowledges the support of BAGEP Award of the Science Academy, Turkey.

\section{Preliminaries}
\label{sec}

Throughout this work $K$ is a field with $\Char K=0$ whose algebraic closure is $\overline{K}$. The Jacobian variety of a smooth algebraic curve $C$ will be denoted by $\Jac(C)$. If two abelian varieties $A$ and $B$ over $K$ are isogenous, we write $A\sim B$. 

 An abelian variety $A$ defined over $K$ is called {\em simple} if there are no lower dimensional abelian varieties $B$ and $C$ over $K$ such that $A$ is isogenous to the product $B \times C$, otherwise it is called {\em decomposable}. If $A$ is simple over $\overline{K}$, then it is called {\em absolutely simple}.

In this note, abusing notation, elliptic curves will be called hyperelliptic curves (of genus $1$). Two hyperelliptic curves of genus $g\ge 2$ described by the following equations
$$y^2 = f(x)\in K[x] \hspace{0.2 cm}\text{and} \hspace{0.3cm} y^2 = f'(x)\in K[x] \hspace{0.2cm}$$
are {\em isomorphic} if and only if
$$x=\dfrac{ax+b}{cx+d}, \quad y=\dfrac{ey}{(cx+d)^{g+1}} , \qquad \hspace{0.2cm} A=\begin{pmatrix}
a & b \\ 
c & d
\end{pmatrix} \in GL_2 (K),\quad \hspace{0.2cm} e \in K^*.$$

  Given a hyperelliptic curve, one would like to know whether its Jacobian is simple or not.

If $A$ is an abelian variety defined over $K$, we write $\End(A)$ for the ring of $\overline{K}$-endomorphisms of $A$. The following results of Zarhin introduce simplicity criteria for certain hyperelliptic Jacobian varieties based on the Galois group of the defining polynomial. 
\begin{Theorem} 
\label{Zar}
Let $C$ be a hyperelliptic curve defined by the equation $y^2= f(x)$, where $f(x)$ is polynomial of degree $ n $ without multiple roots in $K[x]$. 
\begin{itemize}
\item[i)] Assume $n\ge 5$. 
If $\Gal(f)$ is either the full symmetric group $S_n$ or the alternating group $A_n$, then $\End(\Jac(C)) = \Z$. In particular, $\Jac(C)$ is an absolutely simple abelian variety, see \cite{Zarhin1}.
\item[ii)]
Assume $n  \ge 6$ is even. If $f(x) = (x - t)h(x)$ with $t \in K$ and $h(x) \in K[x]$, is such that $\Gal(h)$ is either $S_{n-1}$ or $A_{n-1}$, then $\End(\Jac(C)) = \Z$. In
particular, $\Jac(C)$ is an absolutely simple abelian variety, see \cite{Zarhin2}.
\item[iii)]
Assume $n  \ge 9$ is odd. If $f(x) = (x - t)h(x)$ with $t \in K$ and $h(x) \in K[x]$, is such that $\Gal(h)$ is either $S_{n-1}$ or $A_{n-1}$, then $\Jac(C)$ is an absolutely simple abelian variety, see \cite{Zarhin2}.
\end{itemize}
\end{Theorem}
The following result, \cite[Theorem 8]{Ellenberg} introduces a method to construct absolutely simple varieties over number fields.

\begin{Proposition}
    \label{prop:ellenberg}
    Let $K$ be a
number field. Let $g\ge 1$ be an integer, and let $f\in  K[x]$ be a polynomial of degree $2g$ with no multiple roots.
Consider the hyperelliptic curve of genus $g$ over $K(T)$ defined by $C_T:y^2 = f(x)(x - T)$.
Then there are only finitely many $t \in K$ such that the Jacobian of
$C_t$ is not absolutely simple.
\end{Proposition}

\section{Decomposition into two abelian subvarieties}
Let $C$ be a hyperelliptic curve over $K$ with hyperelliptic involution $\iota$ giving rise to the morphism $C\rightarrow C/\langle \iota\rangle\isom \PP^1 $. We assume that $C$ possesses an automorphism $\sigma$ of order $2$ such that $\sigma\ne \iota$. We set $\tau=\sigma\circ\iota$. Writing $C_{\sigma}$ and $C_{\tau}$ for $C/\langle\sigma\rangle$ and $C/\langle\tau\rangle$ respectively, we obtain the quotient morphisms $\phi_{\sigma}:C\to C_{\sigma}$ and $\phi_{\tau}: C\to C_{\tau}$ respectively. This yields a morphism $\phi=(\phi_{\sigma},\phi_{\tau}):C\to C_{\sigma}\times C_{\tau}$, hence a morphism $\Jac(C)\to \Jac (C_{\sigma})\times \Jac(C_{\tau})$. This morphism is an isogeny, \cite{Kani}, in fact, it is a decomposed Richelot isogeny.

\begin{Lemma} \cite[Theorem 1]{Kat2}
\label{lem:kat}
 Let $C$ be a hyperelliptic curve with an automorphism $\sigma$ of order $2$, which
is not the hyperelliptic involution. We set $\tau = \sigma\circ\iota$ where $\iota$ is the hyperelliptic involution on $C$. Then, the isogeny  $\Jac(C) \rightarrow \Jac(C_{\sigma})\times \Jac(C_{\tau})$ is a decomposed Richelot isogeny.
\end{Lemma}

In this work, we give special attention to the hyperelliptic curve defined by $y^2=f(x^2)$ where $f(x)\in K[x]$ has no multiple roots. 
 
\begin{Proposition}
\label{prop:aut-even}
Let $f(x)\in K[x]\setminus x K[x]$ have no multiple roots. Define the following hyperelliptic curves over $K$
\[C_f:y^2=f(x^2),\qquad E_f:y^2=f(x),\qquad H_f:y^2=xf(x).\]
Then $\Jac(C_f)\sim \Jac(E_f)\times\Jac(H_f)$.
\end{Proposition}
\begin{Proof}
We write $\sigma$ for the automorphism $(x,y)\mapsto (-x,y)$ on $C_f$. The automorphism $\sigma$ is of order $2$. The map $\phi_{\sigma}: C_f\to E_f$ defined by $\phi_{\sigma}:(x,y)\mapsto (x^2,y)$ is the quotient map $C_f\to C_f/\langle\sigma\rangle\isom E_f$. Similarly, if we set $\tau=\sigma\circ\iota$, then $\phi_{\tau}: C_f\to H_f$ defined by $\phi_{\tau}:(x,y)\mapsto (x^2,xy)$ is the quotient map $C_f\to C_f/\langle\tau\rangle\isom H_f$.
\end{Proof}

For an abelian variety $A$ defined over $K$, we set 
$$\End^0(A):=\End(A)\otimes\Q$$ to be the corresponding endomorphism algebra of $A$, which is a semisimple algebra over the field of rational numbers $\Q$. 

\begin{Theorem}
\label{thm:ref}
Let $K$ be a field of characteristic $0$ with algebraic closure $\overline{K}$. Let $f(x) \in K[x]$ be an irreducible polynomial of even
degree $n \ge 8$ such that its Galois group is either the full symmetric
group $S_n$ or the alternating group $A_n$.
Consider the hyperelliptic curve of genus $n-1$ defined by the equation $y^2=f(x^2)$ over $K$. 

Then $\Jac(C_f )$ is isogenous to a product of absolutely simple Jacobian
varietries, $A$ and $B$, of hyperelliptic curves of genus $(n/2-1)$ and $n/2$ respectively. In addition, $\End(A)=\Z$ and $\End(B)=\Z$. In particular, $\End^0(\Jac(C_f))$ is isomorphic to $\Q\oplus\Q$.
\end{Theorem}
\begin{Proof}
Since $f(x)$ is irreducible and $\deg(f) > 1$, one has $f(0) \ne 0$. By
Proposition \ref{prop:aut-even}, there is an isogeny of abelian varieties
$\Jac(C_f ) \sim \Jac (E_f ) \times \Jac(H_f )$ where $E_f$ and $H_f$ are defined as in Proposition \ref{prop:aut-even}. Moreover, the abelian varieties $\Jac(E_f)$ and $\Jac(H_f)$ are of (distinct) dimensions $(n/2-1)$ and $n/2$, respectively. In particular, $\Jac(E_f )$ and
$\Jac(H_f )$ are not isogenous. By Theorem \ref{Zar}, $\Jac(E_f )$ and $\Jac(H_f )$ are absolutely simple. In addition, both endomorphism rings $\End(\Jac(E_f )$
and $\End(\Jac(H_f )$ are the ring of integers $\Z$. Since $\Jac(E_f )$ and $\Jac(H_f )$ are not isogenous, it follows that
the endomorphism algebra of the product $\Jac(E_f )\times \Jac(H_f )$ is isomorphic to  $\Q \oplus \Q$. Now,
since the abelian varieties $\Jac(C_f )$ and $\Jac(E_f ) \times \Jac(H_f )$ are isogenous, their endomorphism algebras are isomorphic. Therefore, the endomorphism algebra of $\Jac(C_f )$ is also isomorphic to $\Q \oplus \Q$. 
\end{Proof}

\begin{Proposition}
\label{prop1}
Let $f(x)\in K[x]$ be of degree $n$ such that $\Gal_K (f)= S_n$ or $A_n$. Let $C_f$, $E_f$ and $H_f$ be as in Proposition \ref{prop:aut-even}.

If $n=2g+1\ge 5$, then $\Jac(C_f)\sim \Jac(E_f)\times\Jac(H_f)$, where both $\Jac(E_f)$ and $\Jac(H_f)$ are absolutely simple of dimension $g$. 

If $n=2g+2\ge 8$, then $\Jac(C_f)\sim \Jac(E_f)\times\Jac(H_f)$, where both $\Jac(E_f)$ and $\Jac(H_f)$ are absolutely simple of dimension $g$ and $g+1$, respectively. 
\end{Proposition}
\begin{Proof}
The statement follows from Proposition \ref{prop:aut-even} and Theorem \ref{Zar}.
\end{Proof}

\begin{Theorem}
\label{thm1}
Let $K$ be a number field. Given any integer $n\ge 2$, there exist infinitely many hyperelliptic curves of genus $n$ with Jacobian varieties that are isogenous over $K$ to the product of two absolutely simple Jacobian varieties of hyperelliptic curves of genus $n/2$ and $n/2$ if $n$ is even; and $(n-1)/2$ and $(n+1)/2$ if $n$ is odd. 
\end{Theorem}
\begin{Proof}
The statement holds in view of Proposition \ref{prop1} for any integer $n$ except possibly $2,3$ and $5$. A hyperelliptic curve with genus two whose Jacobian splits can be constructed easily using Proposition \ref{prop:aut-even}. For example, one may consider the curve $y^2=f(x^2)$ where $f(x)\in K[x]\setminus xK[x]$ is a polynomial of degree $3$ with no multiple roots. 


Let $f(x)$ be a polynomial of degree $d=4$; or of degree $d=6$ with Galois group either $A_6$ or $S_6$. The Jacobian of the curve $y^2=f(x)$ is absolutely simple. This is justified by the fact that the Jacobian is an elliptic curve when $d=4$; or it is an absolutely simple Jacobian of a genus two curve when $d=6$, see Theorem \ref{Zar}. Now, for all but finitely many $t\in K$, the Jacobian of the  curve $y^2=(x-t)f(x)$ is absolutely simple, see Proposition \ref{prop:ellenberg}. For each such value of $t$ such that $t$ is not a root of $f$, we consider the curves $y^2=g_t(x)=f(x+t)$ and $y^2=xg_t(x)$. The latter curves are of genus $1$ and $2$, respectively, when $d=4$; or of genus $2$ and $3$, respectively, when $d=6$, with absolutely simple Jacobians. In addition, the Jacobian of the curve $y^2=g_t(x^2)=f(x^2+t)$ is of dimension $3$ when $d=4$; or of dimension $5$ when $d=6$, for any such $K$-rational value $t$; and it enjoys the required splitting property, see Proposition \ref{prop:aut-even}. 
\end{Proof}

The following proposition indicates that given a polynomial in $K[x]$ of degree $n$ with no multiple roots, one may construct an infinite sequence of hyperelliptic curves of any genus $\ge n-1$ whose Jacobian varieties decompose into two hyperelliptic Jacobian varieties whose dimensions differ by at most $1$.
\begin{Proposition}
    Let $f(x)\in K[x]\setminus xK[x]$ be a polynomial with no multiple roots. Define the following sequence of polynomials 
    \begin{eqnarray*}
    f_0(x)&=& f(x),\\
    g_i(x)&=&x f_i(x),\, i\ge 0,\\
     f_{i}(x)&=& g_{i-1}(x+a_{i-1}),\qquad a_{i-1}\textrm{ is not a root of $g_{i-1}(x)$},\,i\ge 1.   
    \end{eqnarray*}
    Setting $H_{-1}:y^2=f(x)$, $H_i:y^2=g_i(x)$ and $C_i:y^2= f_i(x^2)$, one has
    $\Jac(C_i)\sim \Jac(H_{i-1})\times \Jac(H_i)$, $i\ge 0$. 
    
    If $\deg f=2g+1$, then $H_{i-1}$, $H_i$ and $C_i$ are of genus $g+i/2$, $g+i/2$ and $2g+i$, respectively, when $i$ is even; and of genus $g+r$, $g+r+1$, $2g+i$, respectively, when $i=2r+1$ is odd. 
    
    If $\deg f=2g+2$, then $H_{i-1}$, $H_i$ and $C_i$ are of genus $g+i/2$, $g+i/2+1$, $2g+i+1$, respectively, when $i$ is even; and of genus $g+r+1$, $g+r+1$, $2g+i+1$, respectively, when $i=2r+1$ is odd. 
\end{Proposition}
\begin{Proof}
Observing that $E_i:y^2=f_i(x)$ and $H_{i-1}$, $i\ge 1$, are isomorphic hyperelliptic curves, the proof follows directly from Proposition \ref{prop:aut-even}. 
\end{Proof}
In a similar fashion, we note that the construction of the genus $3$ and $5$ curves using Proposition \ref{prop:ellenberg} in the proof of Theorem \ref{thm1} can be used to provide an alternative way of constructing families of hyperelliptic curves of genus $2n+1\ge 5$ whose Jacobians decompose into the product of two absolutely simple abelian varieties of dimensions $n$ and $n+1$. In addition, the defining polynomials of these curves are essentially multiples of a fixed polynomial of even degree with no multiple roots.  

Given a polynomial $f\in K[x]$ of even degree with no multiple roots, we set $$S(f)=\{t\in K: \textrm{the Jacobian of }y^2=(x-t)f(x) \textrm{ is not absolutely simple; or $t$ is a root of $f(x)$}\}.$$ By Proposition \ref{prop:ellenberg}, $S(f)$ is finite.

\begin{Corollary}
    Let $K$ be a number field. Let $f(x)\in K[x]\setminus xK[x]$ be a polynomial of degree $2g$, $g\ge 1$, with no multiple roots. Define the following sequence of polynomials 
    {\footnotesize{\[\begin{array}{ll}
     f_0(x):=f(x),& g_{0,t_0}(x):= x f_0(x+t_0),\, t_0\not\in S(f_0),\\
    f_{i,t_{i-1}}(x):= (x+r'_{i,t_{i-1}})^{2g+2}g_{i-1,t_{i-1}}\left(\frac{x+r_{i,t_{i-1}}}{x+r'_{i,t_{i-1}}}\right),\, r_{i,t_{i-1}}\ne r'_{i,t_{i-1}},& g_{i,t_i}(x):= x f_{i,t_{i-1}}(x+t_i),\, t_i\not\in S(f_{i,t_{i-1}}), \,i\ge 1,
    \end{array}\]}}
    where $r_{i,t_{i-1}}$ and $r'_{i,t_{i-1}}$ are chosen so that $f_{i,t_{i-1}}(x)\in K[x]\setminus xK[x]$.
    
    Setting $H_{i,t_{i}}:y^2=g_{i,t_i}(x)$, and $C_{i,t_{i-1}}:y^2=f_{i,t_{i-1}}(x^2+t_i)$, then $\Jac(C_{i,t_{i-1}})\sim \Jac(H_{i-1,t_{i-1}})\times\Jac( H_{i,t_{i}})$, where $\Jac( H_{i-1,t_{i-1}})$ is absolutely simple for $i\ge 1$. The genus of the curves $H_{i,t_{i}}$ and $C_{i,t_{i-1}}$ are $g+i$ and $2g+2i-1$, respectively.
\end{Corollary}
\begin{Proof}
We remark that the polynomial $f_{i,t_{i-1}}$ is of even degree. The statement holds in view of Proposition \ref{prop:ellenberg} and Proposition \ref{prop:aut-even} as the curves $H_{i-1,t_{i-1}}$ and $E_{i,t_{i-1}}: y^2=f_{i,t_{i-1}}(x)$ are isomorphic hyperelliptic curves. 
\end{Proof}

\section{Families of decomposable hyperelliptic Jacobian varieties}
We recall that $K$ is a field with $\Char K=0$. We start this section with the following definition. 
\begin{Definition}
    A polynomial $f(x)\in K[x]$ is said to be {\em palindromic} if $f(x)=x^d f(1/x)$ where $d=\deg f$, i.e., if $\displaystyle f(x)=\sum_{i=0}^d a_ix^i$, then $a_i=a_{d-i}$ for $0\le i\le d$. 
\end{Definition}
We write $C_2$, $V_4$ and $D_4$ for the cyclic group with $2$ elements, the Klein-$4$ group, and the dihedral group with $8$ elements, respectively.
\begin{Proposition}
\label{prop:Aut}
Let $f(x)\in K[x]$ be an even palindromic polynomial of degree $2g+2$ with no multiple roots. 
\begin{itemize}
\item[i)] If $C:y^2=f(x)$, then $D_4\hookrightarrow \Aut(C)$, when $g$ is even. 
\item[ii)] If $C:y^2=f(x)$, then $C_2\times C_2\times C_2\hookrightarrow \Aut(C)$, when $g$ is odd. 
\item[iii)] If $C':y^2=xf(x)$, then $V_4\hookrightarrow \Aut(C')$.
\end{itemize}
\end{Proposition}
\begin{Proof}
We write $f(x)=a_{2g+2}x^{2g+2}+a_{2g} x^{2g}+\cdots+a_2x^2+a_0$, where $a_{2i}=a_{2g+2-2i}$, $0\le i\le g+1$. 
For i) and ii) apart from the hyperelliptic involution, the curve $C$ has the following automorphisms of order $2$
\[ \sigma: (x,y)\mapsto (-x,y)\qquad \textrm{and}\qquad \tau:(x,y)\mapsto \left(\frac{1}{x}, \frac{y}{x^{g+1}}\right).\] 
We note that $\sigma^2=\tau^2$. Moreover, $(\sigma\circ\tau)^2=\iota$ when $g$ is even. It follows that the group generated by $\sigma$ and $\tau$ is isomorphic to the dihedral group $D_4$. Specifically,
if we fix a representation $D_4:=\langle a,b| a^2=b^2=(ab)^4=1\rangle$, then we have the following inclusion 
\begin{eqnarray*}
D_4\hookrightarrow\Aut(C);\qquad a\mapsto \sigma, \quad b\mapsto \tau. 
\end{eqnarray*}
ii) follows in a similar fashion by observing that $\sigma\circ\tau$ is an automorphism of order $2$ when $g$ is odd.

 For iii) one my check that the map \[\sigma:C'\rightarrow C':\qquad (x,y)\mapsto \left(\frac{1}{x},\frac{y}{x^{g+2}}\right)\] is an automorphism of $C'$, see \S\ref{sec}, of order $2$. The automorhisms $\iota$, $\sigma$, $\sigma\circ\iota$, $1$ form a subgroup of $\Aut(C')$ isomorphic to the Klein $4$-group, $V_4$.
\end{Proof}

If $f(x)=a_{2g+2}x^{2g+2}+a_{2g} x^{2g}+\cdots+a_2x^2+a_0\in K[x]$ is an even palindromic polynomial with no multiple roots, we write $f_h(x)=a_{2g+2}x^{g+1}+a_{2g} x^{g}+\cdots+a_2x+a_0$. We notice that $f_h(x)$ is a palindromic polynomial itself. We, moreover, set $F_h(x,y)=a_{2g+2}x^{g+1}+a_{2g} x^{g}y+\cdots+a_2xy^g+a_0y^{g+1}$.

\begin{Theorem}
\label{thm:square}
Let $f(x)=a_{2g+2}x^{2g+2}+a_{2g} x^{2g}+\cdots+a_2x^2+a_0\in K[x]$ be an even palindromic polynomial with no multiple roots. Let $f_h(x)=a_{2g+2}x^{g+1}+a_{2g} x^{g}+\cdots+a_2x+a_0$. Assume, moreover, that $C:y^2=f(x)$ and $E:y^2=f_h(x)$. 
\begin{itemize}
\item[i)] If $g\ge 2$ is even, then $\Jac(C)\sim (\Jac(E))^2$. 

\item[ii)] If $g\ge 3$ is odd, then $\Jac(C)\sim \Jac(E)\times\Jac(G_1)\times\Jac(G_2)$ where 
$G_1:y^2=p(x)$ and $G_2:y^2=xp(x)$, and $p(x)\in K[x]$ is such that $p(x^2)=(x^2-1)F_h\left(x+1, x-1\right)$.
\end{itemize}
\end{Theorem}
\begin{Proof}
One observes that $\Jac(C)\sim \Jac(E)\times\Jac(H)$, where $H$ is defined by $y^2=xf_h(x)$, see Proposition \ref{prop:aut-even}. 

If $g=2k$, then $E$ and $H$ are isomorphic hyperelliptic curves via the transformation $$H\longrightarrow E,\quad (x,y)\mapsto \left(\frac{1}{x},\frac{y}{x^{k+1}}\right),$$
see \S\ref{sec}, hence the result.

If $g=2k+1$, then we consider the map 
\begin{eqnarray*}
H\longrightarrow G,\qquad (x,y)\mapsto \left(\frac{x+1}{x-1}, \frac{y}{(x-1)^{k+2}}\right)
\end{eqnarray*}
where $G:y^2=\ell(x)$. One obtains that 
\[\ell(x)=(x^2-1)\left( a_{2g+2}(x+1)^{2k+2}+a_{2g}(x+1)^{2k+1}(x-1)+\cdots+a_2(x+1)(x-1)^{2k+1}+a_0(x-1)^{2k+2}\right),\] hence $\ell(-x)=\ell(x)$, and $\ell$ is an even polynomial of degree $2k+4$. It follows that $\ell(x)=p(x^2)$ for some $p(x)\in K[x]$. 
In view of Proposition \ref{prop:aut-even}, $\Jac(G)\sim \Jac(G_1)\times\Jac(G_2)$, where $G_1:y^2=p(x)$ and $G_2:y^2=xp(x)$. 
\end{Proof}
\begin{Remark}
In Proposition \ref{prop:Aut}, The curve $C':y^2=xf(x)$ possesses the automorphisms $\sigma$ and $\sigma\circ\iota$ described by $(x,y)\mapsto \left(\frac{1}{x},\frac{\pm y}{x^{g+2}}\right)$. In Theorem \ref{thm:square}, the curve $C'$ is described using a different equation, namely, $y^2=p(x^2)$ where the two aforementioned automorphisms are now $(x,y)\mapsto (-x,\pm y)$. Therefore, $C'/\langle\sigma\rangle$ and $C'/\langle\sigma\circ\iota\rangle$ are isomorphic to the hyperelliptic curves defined by $y^2=p(x)$ and $y^2=xp(x)$. 

\end{Remark}

\begin{Corollary}
\label{cor}
\begin{itemize}
    \item[i)] For any integer $n\ge 1$, there exist hyperelliptic curves of genus $2n$ whose Jacobian varieties are isogenous over $K$ to the square of the Jacobian of a hyperelliptic curve of genus $n$. 
    \item[ii)] For any integer $n\ge 1$, there exist hyperelliptic curves of genus $2n+1$ whose Jacobian varieties are isogenous over $K$ to the product of three Jacobian varieties of hyperelliptic curves of genus $n$, $ (n+1)/2$, and $ (n+1)/2$ if $n$ is odd; and $n$, $ 1+n/2$, and $ n/2$ if $n$ is even. 
    \end{itemize}
\end{Corollary}
\begin{Remark}
We remark that Proposition \ref{prop:aut-even} may be used to construct hyperelliptic Jacobian varieties of dimension $2n+1$ that decompose into three Jacobian varieties of lower dimensions, namely, $n+1$, $(n+1)/2$, $(n-1)/2$ if $n$ is odd; and $n+1$, $n/2$, $n/2$ if $n$ is even; which differs from the partitions of the dimension given in Corollary \ref{cor}. In addition, Proposition \ref{prop:aut-even} does not provide a decomposable Jacobian variety whose dimension is $3$.  
\end{Remark}

\begin{Example}
If we consider the curve \[C:y^2=ax^6+bx^4+bx^2+a\in K[x],\] then $\Jac(C)\sim E^2$ where $E$ is the elliptic curve  $y^2=ax^3+bx^2+bx+a.$ 
\end{Example}

\begin{Example}
In Theorem \ref{thm:square}, if one considers the curve 
\[C:y^2=ax^8+bx^6+cx^4+bx^2+a\in K[x]\]
of genus $3$, then $\Jac(C)$ is isogenous to the product of three elliptic curves that are the Jacobians of the following genus $1$ curves
\begin{eqnarray*}
E_1:y^2&=& ax^4+bx^3+cx^2+bx+a,\\
E_2: y^2&=&  (2 a + 2 b + c) x^3+ (10 a - 2 b - 
    3 c) x^2+(-10 a - 2 b + 3 c) x+(-2 a + 2 b - c),\\
E_3: y^2&=&x\left((2 a + 2 b + c) x^3+ (10 a - 2 b - 
    3 c) x^2+(-10 a - 2 b + 3 c) x+(-2 a + 2 b - c)\right) .
\end{eqnarray*}
\end{Example}

\begin{Proposition}
\label{prop4factors}
Let $ f(x)\in K[x]$ be a palindromic polynomial of degree at least $3$. Consider the hyperelliptic  curve $C:y^2=f(x^4)$.  Then $\Jac(C)\sim \Jac(E_1)\times \Jac(E_2) \times \Jac(G_1)\times\Jac(G_2)$, where $E_1:y^2=f(x)$, $E_2:y^2=xf(x)$, and $G_1$ and $G_2$ are as in Theorem \ref{thm:square}.
\end{Proposition}
\begin{Proof}
In view of Theorem \ref{thm:square}, one has $\Jac(C)\sim\Jac(E)\times\Jac(G_1)\times\Jac(G_2)$ where $E:y^2=f(x^2)$. Now due to Proposition \ref{prop:aut-even},  one obtains that $\Jac(E)\sim\Jac(E_1)\times\Jac(E_2)$.
\end{Proof}

\begin{Corollary}
Given any integer $n\ge 2$, there exist hyperelliptic curves of genus $2n+1$ whose Jacobian varieties are isogenous over $K$ to the product of four Jacobian varieties of hyperelliptic curves of genus $(n-1)/2$, $(n+1)/2$, $ (n+1)/2$, and $ (n+1)/2$ if $n$ is odd; and $n/2$, $n/2$, $ n/2$, and $ 1+n/2$ if $n$ is even.  
\end{Corollary}

\begin{Example}
The Jacobian of the hyperelliptic curve $y^2=ax^{12}+bx^8+bx^4+a$ is isogenous to the product of the elliptic curves that are Jacobians of the genus one curves $E_1$, $E_2$, $G_1$; and the Jacobian of the genus $2$ curve $G_2$ 
\begin{eqnarray*}
E_1:y^2&=& ax^3+bx^2+bx+a,\\
E_2:y^2&=& x(ax^3+bx^2+bx+a),\\
G_1:y^2&=&  2 (a + b) x^4+ 2 (14 a - 2 b) x^3+2 (-14 a + 2 b) x+2 (-a - b), \\
G_2:y^2&=& x\left( 2 (a + b) x^4+ 2 (14 a - 2 b) x^3+2 (-14 a + 2 b) x+2 (-a - b)\right).
\end{eqnarray*}
\end{Example}

\section{Rational points on quadratic twists}
In this section, given any integer $g\ge 1$, we construct a hyperelliptic curve of genus $g$ with infinitely many quadratic twists containing at least two $K$-rational non-Weierstrass points. 
\begin{Proposition}
\label{prop:rational1}
Let $f(x)=a_{2g+2}x^{g+1}+a_{2g} x^{g}+\cdots+a_2x+a_0\in K[x]$ be a palindromic polynomial with no multiple roots. Consider the curve $C:y^2=f(x)$. If $g$ is even, then there exists infinitely many quadratic twists of $C$ with at least two $K$-rational non-Weierstrass points.
\end{Proposition}
\begin{Proof}
    Consider the curve $C_{f(t^2)}$ defined over $K(t)$ by
    \[f(t^2)y^2=f(x).\] 
    The set of rational points of $C_{f(t^2)}$ contains the $K(t)$-rational points $\left(t^2,1\right)$ and $\left(\frac{1}{t^2},\frac{1}{t^{2k+1}}\right)$ where $g=2k$. We remark that these points are obtained by considering the quotient maps in \ref{thm:square} i).
\end{Proof}
In the previous proposition, if $g=2$, then $C$ is a genus $1$ curve. Over a number field $K$, this implies the existence of infinitely many quadratic twists of $C$ that are elliptic curves with Mordell-Weil rank at least $2$. That the points are of infinite order follow from Silverman Specialization Theorem, \cite[Theorem 20.3]{Silverman}, whereas the independence of the points follow from the fact that the quotient maps in Theorem \ref{thm:square} are independent maps by construction.

In what follows, we concern ourselves with the construction of tuples of hyperelliptic curves $C_1,\cdots,C_n$ and infinitely many square-free $K$-rational $d$ such that the quadratic twists of these curves by each $d$ contain $K$-rational non-Weierstrass points. 
\begin{Proposition}
\label{prop:rational2}
Let $ f(x)\in K[x]$ be a palindromic polynomial of degree at least $3$ with no multiple roots. Consider the curves $E_1:y^2=f(x)$, $E_2:y^2=xf(x)$, $G_1:y^2=p(x)$ and $G_2:y^2=xp(x)$, where $p(x)$ is defined as in Theorem \ref{thm:square}. There exists infinitely many nonzero $d\in K\setminus K^2$ such that the quadratic twists of $E_1$, $E_2$, $G_1$ and $G_2$ by $d$ contain $K$-rational non-Weierstrass points.  
\end{Proposition}
\begin{Proof} 
We set $n:=\deg f$. We will list down the quadratic twists together with the $K$-rational points on them 
\begin{eqnarray*}
f(t^4)y^2&=&f(x),\qquad (t^4,1),\\
f(t^4)y^2&=& xf(x),\qquad \left(\frac{1}{t^4},\frac{1}{t^{2n+2}}\right),\\
f(t^4)y^2&=& p(x),\qquad \left(\frac{(t^2+1)^2}{(t^2-1)^2},\frac{2^{n+1}t}{(t^2-1)^{n+1}}\right),\\
f(t^4)y^2&=& xp(x),\qquad \left(\frac{(t^2+1)^2}{(t^2-1)^2},\frac{2^{n+1}t(t^2+1)}{(t^2-1)^{n+2}}\right). 
\end{eqnarray*}
These $K$-rational points are obtained using the quotient maps in 
Proposition \ref{prop4factors}.
\end{Proof}
In Proposition \ref{prop:rational2}, if $f$ is chosen to be of degree $3$, then the proposition presents an example of three elliptic curves together with a genus $2$ curve such that there are infinitely many $d$ for which the quadratic twists of these curves by such a $d$ has at least one $K$-rational point. Moreover, if $K$ is a number field, these rational points are of infinite order on the quadratic twists of the elliptic curves, and it is a $K$-rational non-Weierstrass point on the genus two curve.

\end{document}